\documentclass[reqno,centertags, 11pt]{amsart}
\usepackage{amsmath,amsthm,amscd,amssymb}
\usepackage{latexsym}
\newcommand{\bbC}{{\mathbb{C}}}
\newcommand{\bbD}{{\mathbb{D}}}
\newcommand{\bbE}{{\mathbb{E}}}

\newcommand{\bbR}{{\mathbb{R}}}
\newcommand{\bbT}{{\mathbb{T}}}

\newcommand{\calE}{{\mathcal E}}

\newcommand{\dott}{\,\cdot\,}

\newcommand{\lb}{\label}
\newcommand{\f}{\frac}

\newcommand{\ol}{\overline}
\newcommand{\ti}{\tilde  }

\newcommand{\ac}{\text{\rm{ac}}}
\newcommand{\s}{\text{\rm{s}}}

\newcommand{\bi}{\bibitem}

\newcommand{\beq}{\begin{equation}}
\newcommand{\eeq}{\end{equation}}
\newcommand{\ba}{\begin{align}}
\newcommand{\ea}{\end{align}}
\newcommand{\veps}{\varepsilon}

\newcounter{smalllist}
\newenvironment{SL}{\begin{list}{{\rm\roman{smalllist})}}{%
\setlength{\topsep}{0mm}\setlength{\parsep}{0mm}\setlength{\itemsep}{0mm}%
\setlength{\labelwidth}{2em}\setlength{\leftmargin}{2em}\usecounter{smalllist}%
}}{\end{list}}

\newcommand{\bigtimes}{\mathop{\mathchoice%
{\smash{\vcenter{\hbox{\LARGE$\times$}}}\vphantom{\prod}}%
{\smash{\vcenter{\hbox{\Large$\times$}}}\vphantom{\prod}}%
{\times}%
{\times}%
}\displaylimits}

\DeclareMathOperator{\Real}{Re}
\DeclareMathOperator{\Ima}{Im}
\allowdisplaybreaks
\numberwithin{equation}{section}

\newtheorem{theorem}{Theorem}[section]
\newtheorem{proposition}[theorem]{Proposition}

\theoremstyle{definition}

\theoremstyle{remark}

\newcommand{\abs}[1]{\lvert#1\rvert}

\begin{document}
\title[Analogs of the M-Function in OPUC Theory]
{Analogs of the M-Function in the Theory of Orthogonal Polynomials on the Unit Circle}
\author[B. Simon]{Barry Simon$^*$}

\thanks{$^*$ Mathematics 253-37, California Institute of Technology, Pasadena, CA 91125, USA. 
E-mail: bsimon@caltech.edu. Supported in part by NSF grant DMS-0140592} 
\thanks{{\it For special issue of J. Comp. Appl. Math.}}

\dedicatory{To Norrie Everitt, on his 80th birthday, \\ 
a bouquet to the master of the $m$-function}

\date{October 1, 2003}

\begin{abstract} We show that the multitude of applications of the Weyl-Titchmarsh $m$-function 
leads to a multitude of different functions in the theory of orthogonal polynomials on the unit 
circle that serve as analogs of the $m$-function. 
\end{abstract}

\maketitle

\section{Introduction} \lb{s1}

Use of the Weyl-Titchmarsh $m$-function has been a constant theme in Norrie Everitt's 
opus, so I decided a discussion of the analogs of these ideas in the theory of 
orthogonal polynomials on the unit circle (OPUC) was appropriate. Interestingly enough, 
the uses of the $m$-functions are so numerous that OPUC has multiple analogs 
of the $m$-function! 

$m$-functions are associated to solutions of 
\begin{equation} \lb{1.1} 
-u'' + qu =zu
\end{equation}
with $q$ a real function on $[0,\infty)$ and $z$ a parameter in $\bbC_+ =\{z\mid 
\Ima z >0\}$. The most fundamental aspect of the $m$-function is its relation to the 
spectral measure, $\rho$, for \eqref{1.1} by 
\begin{equation} \lb{1.2}
m(z)=c+\int d\rho(x) \biggl[\f{1}{x-z} - \f{x}{1+x^2}\biggr]
\end{equation}
where $c$ is determined by (see Atkinson \cite{Atk81}, Gesztesy-Simon \cite{GS}):
\begin{equation} \lb{1.3}
m(z) =\sqrt{-z} + o(1) \qquad \text{as } z\to i\infty
\end{equation}

\eqref{1.2} plus \eqref{1.3} allow you to compute $m$ given $d\rho$, and $d\rho$ is 
determined by $m$ via 
\begin{equation} \lb{1.4}
\lim_{e\downarrow 0}\, \f{1}{\pi} \int_a^b m(x+i\veps)\, dx = \tfrac12 [\rho((a,b)) 
+ \rho([a,b])]
\end{equation}

Of course, I haven't told you what $m$ or $\rho$ is. This is done by defining $m$, in 
which case $\rho$ is defined by \eqref{1.4}. Under weak conditions on $q$ at $\infty$, 
for $z\in\bbC_+$, \eqref{1.1} has a solution $u(x,z)$ which is $L^2$ at infinity, and 
it is unique up to a constant multiple. Then, $m$ is defined by 
\begin{equation} \lb{1.5}
m(z) =\f{u'(0,z)}{u(0,z)}
\end{equation}
With this definition, $d\rho$ is a spectral measure for $u\mapsto -u'' + qu =Hu$ in the 
sense that $H$ is unitarily equivalent to multiplication by $\lambda$ on $L^2 (\bbR,d\rho)$.  
\eqref{1.5} is often written in the equivalent form,  
\[
\psi (x,z) +m(z) \varphi(x,z)\in L^2 
\]
where $\varphi,\psi$ solve \eqref{1.1} with initial conditions $\varphi(0)=0$, $\varphi' 
(0)=1$, $\psi(0)=1$, $\psi'(0)=0$. 

Note that if one defines 
\begin{equation} \lb{1.6}
m(x;z) = \f{u'(x,z)}{u(x,z)} 
\end{equation}
the $m$-function for $q_x (\cdot) =q(\dott + x)$, then $m$ obeys the Riccati equation 
\begin{equation} \lb{1.7}
m' = q-z-m^2
\end{equation} 

It could be said that this is backwards: the definition \eqref{1.5} should come first, 
before \eqref{1.2}. I put it in this order because it is \eqref{1.2} that makes $m$ such 
an important object both in classical results \cite{Atk64,CL,East73,East02,EKZ,Hart,Reid,Ti} 
and very recent work \cite{S271,Gel,Lev,Stie,SFI,Borg}. 

To describe the third role of the $m$-function, it will pay to switch to the case of 
Jacobi matrices. We now have, instead of $q$, two sequences $\{a_n\}_{n=1}^\infty$, 
$\{b_n\}_{n=1}^\infty$ with $a_n >0$, $b_n\in\bbR$ which we will suppose uniformly 
bounded. Define an infinite matrix 
\begin{equation} \lb{1.8x}
J= \begin{pmatrix}
b_1 & a_1 & 0 & 0 & \cdots \\
a_1 & b_2 & a_2 & 0 & \cdots \\
0 & a_2 & b_3 & a_3 & \cdots \\
\vdots & \vdots & \vdots & \vdots & \ddots
\end{pmatrix}
\end{equation} 
which is a bounded selfadjoint operator. One defines 
\begin{equation} \lb{1.8}
m(z) =\langle \delta_1, (J-z)^{-1}\delta_1\rangle 
\end{equation} 
In terms of the spectral measure, $\mu$, for $\delta_1$ for $J$, 
\begin{equation} \lb{1.8a}
m(z) =\int \f{d\mu(x)}{x-z}
\end{equation} 
If $u_n$ is the $\ell^2$ solution of $a_{n-1} u_{n-1} + (b_n -z) u_n + a_n u_{n+1} =0$ 
with $\Ima z >0$, one has the analog of \eqref{1.5} 
\begin{equation} \lb{1.8e}
m(z) =\f{u_1(z)}{u_0(z)}
\end{equation}

This process of going from $a$ and $b$ to $m$ and then to $\mu$ can be reversed. 
One way is by iterating \eqref{1.10} below, which lets one go from $\mu$ to $m$ 
(by \eqref{1.8a}) and then gets the $a$'s and $b$'s as coefficients in a continued 
fraction expansion of $m$. From our point of view, an even more important way of 
going backwards uses orthogonal polynomials on the real line (OPRL). Given $\mu$ 
(of bounded support), one forms the monic orthogonal polynomials $P_n(x)$ for $d\mu$ 
and shows they obey a recursion relation 
\begin{equation} \lb{1.8b}
P_{n+1}(x) = (x-b_{n+1}) P_n(x) -a_n^2 P_{n-1}(x)
\end{equation}
which yields the Jacobi parameters $a$ and $b$. The orthonormal polynomials, $p_n(x)$, 
are related to $P_n$ by 
\begin{equation} \lb{1.8c}
p_n(x) = (a_1 \dots a_n)^{-1} P_n(x)
\end{equation}
and obey 
\begin{equation} \lb{1.8d}
a_{n+1} p_{n+1}(x) = (x-b_{n+1}) p_n(x) -a_n p_{n-1}(x)
\end{equation} 

\eqref{1.7} has the analog 
\begin{equation} \lb{1.10}
m(z;J) = (b_1 -z- a_1^2 m(z;J^{(1)}))^{-1}
\end{equation}
where $J^{(1)}$ is the Jacobi matrix with parameters $\ti a_m =a_{m+1} \ti b_m = b_{m+1}$ 
(i.e., the top row and left column are removed). 

If $m(x+i\veps; J)$ has a limit as $\veps\downarrow 0$, \eqref{1.10} says that $m(x+i\veps; 
J^{(1)})$ has a limit, and by \eqref{1.10}, 
\begin{equation} \lb{1.11}
\f{\Ima m(x;J)}{\Ima m(x;J^{(1)})} = \abs{a_1 m(x;J)}^2
\end{equation} 
$\Ima m$ is important because if $\mu$ is given by \eqref{1.8a} then 
\begin{equation} \lb{1.12a}
d\mu_\ac = \f{1}{\pi}\, \Ima m(x+i0)\, dx 
\end{equation}
This property of $m$, that its energy is the ratio of $\Ima$'s, is a critical element of 
recent work on sum rules for spectral theory \cite{OPUC,KS,SZ,Sim288,Jost1}. 

The interesting point is that, for OPUC, the analogs of the functions obeying \eqref{1.2}, 
\eqref{1.5}, and \eqref{1.11} are different! In Section~\ref{s2}, we will give a quick 
summary of OPUC. In Section~\ref{s3}, we discuss \eqref{1.2}; in Section~\ref{s4}, we 
discuss \eqref{1.11}; and finally, in Section~\ref{s5}, the analog of \eqref{1.5}. 

\medskip

Happy 80th, Norrie. I hope you enjoy this bouquet.

\medskip
\section{Overview of OPUC} \lb{s2} 

We want to discuss here the basics of OPUC, although we will only scratch the surface of a 
rich and beautiful subject \cite{OPUC}. The theory reverses the usual passage from 
differential/difference equations to measures, and instead follows the discussion of OPRL 
in Section~\ref{s1}. $\mu$ is now a probability measure on $\partial\bbD =\{z\mid \abs{z}=1\}$. 
We suppose $\mu$ is nontrivial, that is, not supported on a finite set. One can then form, 
by the Gram-Schmidt procedure, the monic orthogonal polynomials $\Phi_n(z)$ and the 
orthonormal polynomials, $\varphi_n(z) = \Phi_n(z)/\|\Phi_n\|$ where $\|\cdot\|$ is the 
$L^2 (\partial\bbD,d\mu)$ norm. 

Given fixed $n\in \{0,1,2,\dots \}$, we define an anti-unitary operator on $L^2 
(\partial\bbD, d\mu)$ by 
\begin{equation} \lb{2.1}
f^*(z) =z^n \, \ol{f(z)}
\end{equation} 
The use of a symbol without ``$n$" is terrible notation, but it is standard! If $Q_n$ 
is a polynomial of degree $n$, $Q_n^*$ is also a polynomial of degree $n$. Indeed, 
\[
Q_n^*(z) = z^n\, \ol{Q_n (1/\bar z)}
\]
so if $Q_n(z) =a_n z^n + a_{n-1} z^{n-1} + \cdots + a_0$, then $Q_n^*(z) = \bar a_0 z^n 
+ \bar a_1 z^{n-1} + \cdots + \bar a_n$. 

Since $\Phi_n$ is monic, $\Phi_n^*(0)=1$, and thus, $N(z) \equiv (\Phi_{n+1}^*(z) - 
\Phi_n^*(z))/z$ is a polynomial of degree $n$. Since $\,{}^*\,$ is anti-unitary, 
\begin{align*} 
\langle z^m, N(z)\rangle &= \langle z^{m+1}, \Phi_{n+1}^* - \Phi_n^*\rangle \\
&= \langle \Phi_{n+1}, z^{n+1 -(m+1)}\rangle - \langle\Phi_n, z^{n-m-1}\rangle \\
&= 0
\end{align*} 
for $m=0,1,\dots, n-1$. Thus $N(z)$ must be a multiple of $\Phi_n(z)$, that is, for some 
$\alpha_n\in\bbC$, 
\begin{equation} \lb{2.2}
\Phi_{n+1}^*(z) = \Phi_n^*(z) -\alpha_n z\Phi_n(z)
\end{equation}
and its $\,{}^*\,$,
\begin{equation} \lb{2.3}
\Phi_{n+1}(z) = z\Phi_n(z) -\bar\alpha_n \Phi_n^*(z)
\end{equation}
\eqref{2.2}/\eqref{2.3} are the {\it Szeg\H{o} recursion formulae} (\cite{Szb}); the 
$\alpha_n$'s are the Verblunsky coefficients (after \cite{V35}). The derivation I've just 
given is that of Atkinson \cite{Atk64}. 

Since $\Phi_n^*\perp\Phi_{n+1}$, \eqref{2.3} implies 
\[
\|\Phi_{n+1}\|^2 + \abs{\alpha_n}^2 \|\Phi_n^*\|^2 = \|z \Phi_n\|^2
\]
Since $\|\Phi_n^*\| =\|z\Phi_n\| = \|\Phi_n\|$, we have 
\begin{equation} \lb{2.4}
\|\Phi_{n+1}\| = (1-\abs{\alpha_n}^2)^{1/2} \|\Phi_n\|
\end{equation}
This implies first of all that 
\begin{equation} \lb{2.5}
\abs{\alpha_n} <1
\end{equation}
and if 
\begin{equation} \lb{2.6}
\rho_n\equiv (1-\abs{\alpha_n}^2)^{1/2}
\end{equation}
then 
\begin{equation} \lb{2.7}
\|\Phi\|_n = \rho_0\rho_1 \dots \rho_{n-1}
\end{equation}
so 
\begin{equation} \lb{2.8}
\varphi_n =(\rho_0 \dots \rho_{n-1})^{-1} \Phi_n
\end{equation} 
and \eqref{2.2}, \eqref{2.3} becomes 
\begin{align} 
z\varphi_n &=\rho_n \varphi_{n+1} + \bar\alpha_n \varphi_n^* \lb{2.9} \\ 
\varphi_n^* &=\rho_n \varphi_{n+1}^* + \alpha_n z\varphi_n \lb{2.10}
\end{align} 

The $\alpha_n$'s not only lie in $\bbD$, but it is a theorem of Verblunsky \cite{V35} 
that as $\mu$ runs through all nontrivial measures, the set of $\alpha$'s runs through 
all of $\bigtimes_{n=0}^\infty \bbD$. The $\alpha$'s are the analogs of the $a$'s and 
$b$'s in the Jacobi case or of $V$ in the Schr\"odinger case.  

We will later have reason to consider Szeg\H{o}'s theorem in Verblunsky's form 
\cite{V36}: 

\begin{theorem}\lb{T2.1} Let 
\begin{equation} \lb{2.10a}
d\mu = w\, \f{d\theta}{2\pi} + d\mu_\s
\end{equation}
Then 
\begin{equation} \lb{2.11}
\prod_{j=0}^\infty (1-\abs{\alpha_j}^2) = \exp \biggl( \int \log (w(\theta)) \f{d\theta}{2\pi}\biggr)
\end{equation}
\end{theorem} 

{\it Remark.} The log integral can diverge to $-\infty$. The theorem says the integral is 
$-\infty$ if and only if the product on the left is $0$, that is, if and only if $\sum 
\abs{\alpha_j}^2 =\infty$. 

If 
\begin{equation} \lb{2.12}
\sum_{j=0}^\infty\, \abs{\alpha_j}^2 <\infty 
\end{equation}
we say the Szeg\H{o} condition holds. This happens if and only if 
\begin{equation}\lb{2.13} 
\int \abs{\log (w(\theta))}\, \f{d\theta}{2\pi} <\infty
\end{equation}
In that case, we define the Szeg\H{o} function on $\bbD$ by 
\begin{equation} \lb{2.14}
D(z) =\exp \biggl(\int \f{e^{i\theta}+z}{e^{i\theta}-z} \, \log (w(\theta)) \, 
\f{d\theta}{4\pi}\biggr)
\end{equation}

\medskip

\section{The Carath\'eodory and Schur Functions} \lb{s3} 

Given \eqref{1.8a} (and \eqref{1.2}), the natural ``$m$-function" for OPUC is the 
Carath\'eodory function, $F(z)$, 
\begin{equation} \lb{3.1}
F(z) =\int \f{e^{i\theta}+z}{e^{i\theta}-z}\, d\mu(\theta)
\end{equation}
The Cauchy kernel $(e^{i\theta}+z)/(e^{i\theta}-z)$ has the Poisson kernel 
\begin{equation} \lb{3.2}
\left. \Real \biggl( \f{e^{i\theta}+z}{e^{i\theta}-z}\biggr)
\right|_{z=re^{i\varphi}} = \f{1-r^2}{1+r^2 - 2\cos (\theta-\varphi)}
\end{equation}
as its real part, and this is positive, so 
\begin{equation} \lb{3.3}
\Real F(z) >0 \text{ for }z\in\bbD \qquad F(0)=1
\end{equation}
This replaces $\Ima m>0$ if $\Ima z>0$. 

One might think the ``correct" analog of $m$ is 
\begin{equation} \lb{3.4}
R(z) =\int \f{1}{e^{i\theta}-z}\, d\mu(\theta)
\end{equation}
$R$ and $F$ are related by 
\begin{equation} \lb{3.5}
R(z) =(2z)^{-1} (F(z)-1)
\end{equation}

If one rotates $d\mu$ and $z$ (i.e., $d\mu(\theta)\to d\mu (\theta-\varphi)$, 
$z\to e^{i\varphi}z$), $F$ is unchanged but $R$ is multiplied by $e^{-i\varphi}$,  
so the set of values $R$ can take are essentially arbitrary --- which shows $F$, which 
obeys $\Real F(z)>0$, is a nicer object to take. That said, we will see $R$ again in 
Section~\ref{s5}. 

$F$ has some important analogs of $m$: 
\begin{SL} 
\item[(1)] $\lim_{r\uparrow 1} F(re^{i\theta})$ exists for a.e.~$\theta$, and if 
\eqref{2.10a} defines $w$, then 
\begin{equation} \lb{3.6} 
w(\theta) =\Real F(e^{i\theta}) 
\end{equation}
\item[(2)] $\theta_0$ is a pure point of $\mu$ if and only if $\lim_{r\uparrow 1} (1-r) 
\Real F(re^{i\theta_0})\neq 0$ and, in general, 
\[
\lim_{r\uparrow 1}\, (1-r) \Real F(re^{i\theta_0}) =\mu(\{\theta_0\})
\]
\item[(3)] $d\mu_\s$ is supported on $\{\theta\mid\lim_{r\uparrow 1} F(re^{i\theta}) 
=\infty\}$.
\end{SL} 
In fact, the proof of the analogs of these facts for $m$ proceeds by mapping $\bbC_+$ 
to $\bbD$ and using these facts for $F$! 

These properties provide a strong analogy, but one can note a loss of ``symmetry" 
relative to the ODE case. The $m$-function maps $\bbC_+$ to $\bbC_+$. $F$ though maps 
$\bbD$ to $-i\bbC_+$. One might prefer a map of $\bbD$ to $\bbD$. In fact, one defines the 
Schur function, $f$, of $\mu$ via 
\begin{equation} \lb{3.7} 
F(z) =\f{1+zf(z)}{1-zf(z)} 
\end{equation}
then $f$ maps $\bbD$ to $\bbD$ and \eqref{3.7} sets up a one-one correspondence between 
$F$'s with $\Real F>0$ on $\bbD$ and $F(0)=1$ and $f$ mapping $\bbD$ to $\bbD$ (this fact 
relies on the Schwarz lemma that $f$ maps $\bbD$ to $\bbD$ with $f(0)=0$ if and only if 
$f=zg$ where $g$ maps $\bbD$ to $\bbD$). 

For at least some purposes, $f$ is a ``better" analog of $m$ than $F$, for example, in 
regard to its analog of the recursion \eqref{1.10}. If $f$ is the Schur function associated 
to Verblunsky coefficients $\{\alpha_0, \alpha_1, \dots\}$ and $f_n$ is the Schur function 
associated to $\{\alpha_n, \alpha_{n+1}, \dots\}$, then 
\begin{equation} \lb{3.8} 
f=\f{\alpha_0 + zf_1}{1+\bar\alpha_0 zf_1}
\end{equation}
a result of Geronimus (see \cite{OPUC} for lots of proofs of this fact). 

Interestingly enough, Schur, not knowing of the connection to OPUC, discussed \eqref{3.8} 
for $\alpha_0 =f(0)$ as a map of $f\to (\alpha_0, f_1)$ and, by iteration, to a 
parametrization of functions of $\bbD$ to $\bbD$ by parameters $\alpha_0, \dots, \alpha_n, 
\dots$. There is, of course, a formula relating $F$ to $F_1$ that can be obtained from 
\eqref{3.7} and \eqref{3.8} or directly \cite{Pe96}, but it is more complicated than 
\eqref{3.8}. 

Finally, in discussing $f$, we note that there is a natural family 
$\{d\mu_\lambda\}_{\lambda\in\partial\bbD}$ of measures related to $d\mu$ (with 
$d\mu_{\lambda =1}=d\mu$) that corresponds to ``varying boundary conditions." We will 
discuss those more fully in Section~\ref{s5}, but we note 
\begin{equation} \lb{3.9} 
f(z;d\mu_\lambda) =\lambda f(z;d\mu) 
\end{equation}
while the formula for $F(d\mu_\lambda)$ is more involved. 

The Schur function and Schur iterates, $f_n$, have been used by Khrushchev 
\cite{Kh2000,Khr,KhGo} as a powerful tool in the analysis of OPUC.  

\medskip
\section{The Relative Szeg\H{o} Function}\lb{s4} 

As explained in the introduction, a critical property of $m$ is \eqref{1.11}, which is the 
basis of step-by-step sum rules (see \cite{Sim288}). The left side of \eqref{1.11} enters 
as the ratio of a.c.~weights of $d\mu_J$ and $d\mu_{J^{(1)}}$. Thus, we are interested in 
$\Ima F(e^{i\theta};\{\alpha_j\}_{j=0}^\infty)$ divided by $\Ima F(e^{i\theta}; 
\{\alpha_{j+1}\}_{j=0}^\infty)$, that is, $\Ima F/\Ima F_1$ in the language of the last 
section. Neither $\abs{F}$ nor $\abs{f}$ is directly related to this ratio, so we need a 
different object to get an analog of \eqref{1.11}. The following was introduced by 
Simon in \cite{OPUC}:
\begin{equation} \lb{4.1} 
(\delta_0D)(z) =\f{1-\bar\alpha_0 f}{\rho_0}\,\, \f{1-zf_1}{1-zf}
\end{equation} 
It is called the ``relative Szeg\H{o} function" for reasons that will become clear in a 
moment. 

In \eqref{4.1}, $f_1$ is the Schur function for Verblunsky coefficients 
\begin{equation} \lb{4.2} 
\alpha_j^{(1)} =\alpha_{j+1} 
\end{equation}

Here is the key fact: 

\begin{theorem}\lb{T4.1} Let $d\mu$ and $d\mu^{(1)}$ be measures on $\partial\bbD$ with 
Verblunsky coefficients related by \eqref{4.2}. Suppose $d\mu =w(\theta) \f{d\theta}{2\pi} + 
d\mu_\s$ and $d\mu^{(1)}=w^{(1)}\f{d\theta}{2\pi} +d\mu_\s$. Then 
\begin{SL} 
\item[{\rm{(1)}}] For a.e.~$\theta$, $\lim_{r\uparrow 1} (\delta_0D)(re^{i\theta}) \equiv 
\delta_0 D (e^{i\theta})$ exists. 
\item[{\rm{(2)}}] If $w(\theta) \neq 0$, then {\rm{(}}for a.e.~$\theta$ w.r.t.~$\f{d\theta}
{2\pi}${\rm{)}}, $w_1(\theta)\neq 0$ and 
\begin{equation} \lb{4.3} 
\f{w(\theta)}{w_1(\theta)} =\abs{(\delta_0 D)(e^{i\theta})}^2
\end{equation}
\end{SL}
\end{theorem} 

\begin{proof}[Sketch of Proof] Each of the functions $1-\bar\alpha_0 f$, $1-zf_1$, and 
$1-zf$ takes values in $\{w\mid \abs{w-1}<1\}$ on $\bbD$, so their arguments lie in 
$[-\f{\pi}{2},\f{\pi}{2}]$, so their logs are in all $H^p$, $1<p<\infty$. That is, they 
are outer functions, and so $\delta_0 D$ is an outer function, which means that 
assertion (1) holds (see Rudin \cite{Rudin} for a pedagogic discussion of outer functions). 

To get \eqref{4.3}, we note that \eqref{3.7} implies 
\[
\Real F(z) = \f{1-\abs{f}^2\abs{z}^2}{\abs{1-zf}^2} 
\]
so 
\begin{equation} \lb{4.4} 
\f{\Real F(z)}{\Real F_1(z)} = \biggl| \f{1-zf_1}{1-zf}\biggr|^2 
\f{1-\abs{f}^2 \abs{z}^2}{1-\abs{f_1}^2 \abs{z}^2} 
\end{equation}

On the other hand, \eqref{3.8} implies 
\begin{equation} \lb{4.5}  
zf_1 = \f{f-\alpha_0}{1-\bar\alpha_0 f} 
\end{equation}
which implies 
\begin{equation} \lb{4.6} 
1-\abs{zf_1}^2 = \f{\rho_0^2 (1-\abs{f}^2)}{\abs{1-\bar\alpha_0 f}^2} 
\end{equation}
so, putting these formulae together, 
\begin{equation} \lb{4.7} 
\f{\Real F(z)}{\Real F_1(z)} = \abs{(\delta_0 D)(z)}^2 
\biggl( \f{1-\abs{z}^2 \abs{f}^2}{1-\abs{f}^2}\biggr)
\end{equation}
which, as $\abs{z}\to 1$, yields \eqref{4.3}. 
\end{proof} 

In particular, one has the nonlocal step-by-step sum rule that if $w(\theta) \neq 0$ for 
a.e.~$\theta$, then 
\begin{equation} \lb{4.8} 
(\delta_0 D)(z) =\exp\biggl( \int_0^{2\pi} \f{e^{i\theta}+z}{e^{i\theta}-z} \, 
\log \biggl( \f{w(\theta)}{w_1(\theta)}\biggr) \, \f{d\theta}{4\pi}\biggr) 
\end{equation}
and, in particular, setting $z=0$, 
\begin{equation} \lb{4.9} 
\rho_0^2 = \exp \biggl( \int_0^{2\pi} \log \biggl( \f{w(\theta)}{w_1(\theta)}\biggr) 
\, \f{d\theta}{2\pi}\biggr) 
\end{equation}
which is not only consistent with Szeg\H{o}'s theorem \eqref{2.11} but, using 
semicontinuity of the entropy, can be used to prove it (see \cite{KS,OPUC}) as 
follows: 
\begin{SL} 
\item[(1)] Iterating \eqref{4.9} yields 
\begin{equation} \lb{4.10} 
(\rho_0 \dots \rho_{n-1})^2 = \exp \biggl( \int_0^{2\pi} \log 
\biggl(\f{w(\theta)}{w_n(\theta)}\biggr)\, \f{d\theta}{2\pi}\biggr) 
\end{equation}
\item[(2)] Since $\exp (\int_0^{2\pi} \log(w_n(\theta)\f{d\theta}{2\pi}) 
\leq \int_0^{2\pi} w_n (\theta) \f{d\theta}{2\pi} \leq 1$, \eqref{4.10} implies 
\begin{equation} \lb{4.11} 
(\rho_0 \dots \rho_{n-1})^2 \geq \exp \biggl( \int_0^{2\pi} \log 
(w(\theta))\, \f{d\theta}{2\pi}\biggr) 
\end{equation}
\item[(3)] If $w^{(n)}$ is the weight associated to the measure with 
\[
\alpha_j^{(n)} =\begin{cases} 
\alpha_j & j\leq n -1 \\
0 & j\geq n \end{cases}
\]
\eqref{4.10} proves 
\begin{equation} \lb{4.12} 
(\rho_0  \dots \rho_{n-1})^2 = \exp \int _0^{2\pi} \log (w^{(n)}(\theta)) \, 
\f{d\theta}{2\pi} 
\end{equation}
\item[(4)] $d\mu\to\int_0^{2\pi} \log (w(\theta))\f{d\theta}{2\pi}$ is an entropy,  
hence, weakly upper semicontinuous. Since $w^{(n)}\f{d\theta}{2\pi} \to d\mu$ 
weakly as $n\to\infty$, this semicontinuity shows 
\begin{equation} \lb{4.13} 
\lim_{n\to\infty} \, (\rho_n \dots \rho_{n-1})^2 \leq \exp \biggl( \int_0^{2\pi} 
\log (w(\theta)) \, \f{d\theta}{2\pi}\biggr) 
\end{equation}
\eqref{4.11} and \eqref{4.13} is Szeg\H{o}'s theorem. 
\end{SL} 

\smallskip
Two other properties of $\delta_0 D$ that we should mention are: 
\begin{SL} 
\item[(A)] If $\sum_{n=0}^\infty \abs{\alpha_n}^2 <\infty$, then 
\begin{equation} \lb{4.14} 
(\delta_0 D)(z) = \f{D(z;\alpha_0, \alpha_1, \alpha_2, \dots)}
{D(z;\alpha_1, \alpha_2, \alpha_3, \dots)} 
\end{equation}
\item[(B)] In general, one has 
\begin{equation} \lb{4.15} 
\delta_0 D(z) = \lim_{n\to\infty} \, \f{\varphi_{n-1}^*(z;\alpha_1, \alpha_2, \dots)} 
{\varphi_n^* (z;\alpha_0, \alpha_1, \dots)}
\end{equation}
\end{SL}

\medskip
\section{Eigenfunction Ratios}\lb{s5} 

Finally, we look at the analogs of $m$ as a function ratio, its initial definition by 
Weyl and Titchmarsh. The key papers on this point of view are by Geronimo-Teplyaev 
\cite{GTep} and Golinskii-Nevai \cite{GN}. We will see from one point of view \cite{GN}  
that $F(z)$ plays this role, but from other points of view \cite{GTep} that other 
functions are more natural. 

The recursion relations \eqref{2.9}/\eqref{2.10} can be rewritten as 
\begin{equation} \lb{5.1} 
\binom{\varphi_{n+1}}{\varphi_{n+1}^*} = A(\alpha_n,z) \binom{\varphi_n}{\varphi_n^*} 
\end{equation}
where 
\begin{equation} \lb{5.2} 
A(\alpha,z) = \rho^{-1} \begin{pmatrix} 
z & -\bar\alpha_n \\ -\alpha_n z & 1 \end{pmatrix}  
\end{equation}
(with $\rho=(1-\abs{\alpha}^2)^{1/2}$). From this point of view, the analog of the 
fundamental differential/difference equation in the real case is 
\begin{equation}\lb{5.3}
\Xi_n =T_n (z) \Xi_0
\end{equation}
with 
\begin{equation}\lb{5.4}
T_n(z) = A(\alpha_{n-1}, z) \dots A(\alpha_0, z)
\end{equation}
The correct boundary conditions for the usual OPUC are $\Xi_0 =\binom{1}{1}$. 

One can ask for what other initial conditions the polynomials associated with 
the top component of $T_n (z) \Xi_0$ are OPUC for some measure. Note that 
\begin{equation}\lb{5.5}
\binom{1}{\lambda} = U(\lambda) \binom{1}{1}
\end{equation}
with 
\begin{equation}\lb{5.6}
U(\lambda) = \begin{pmatrix} 1 & 0 \\ 0 & \lambda \end{pmatrix}
\end{equation}
and that 
\begin{equation}\lb{5.7}
U(\lambda)^{-1} A(\alpha,z) U(\lambda) =\rho^{-1} 
\begin{pmatrix} z & -\bar\alpha_n \lambda \\ -\alpha_n \lambda^{-1} z & 1 
\end{pmatrix}
\end{equation} 

We see from this that $\bar\lambda =\lambda^{-1}$, that is, $\abs{\lambda} =1$ 
will yield $U(\lambda)^{-1} A(\alpha_1, z) U(\lambda) = A(\bar\lambda \alpha, z)$. 
Changing $\lambda$ to $\bar\lambda$, we see that 

\begin{proposition}\lb{P5.1} Let $\abs{\lambda} =1$. If $\varphi_n^{(\lambda)}(z)$ 
are the OPUC for Verblunsky coefficients $\alpha_n^{(\lambda)}=\lambda\alpha_n$, then 
\begin{equation}\lb{5.8}
\binom{\varphi_n^{(\lambda)}(z)}{\bar\lambda \varphi_n^{(\lambda)*}(z)} = 
T_n (z;\{\alpha_j\}_{j=1}^\infty) \binom{1}{\bar\lambda}
\end{equation}
\end{proposition} 

This suggests that one look at the family $d\mu_\lambda$ or measures with 
\begin{equation}\lb{5.9}
\alpha_j (d\mu_\lambda) =\lambda \alpha_j (d\mu)
\end{equation}
called the family of Aleksandrov measures associated to $\{\alpha_j\}_{j=0}^\infty$ 
after \cite{Alex}. The special case $\lambda =-1$ goes back to Verblunsky \cite{V36}  
and Geronimus \cite{Ger46}, and are called the second kind polynomials, denoted 
$\psi_n(z)$. The following goes back to Verblunsky \cite{V36}: 

\begin{theorem}\lb{T5.2} For $z\in\bbD$, uniformly on compact subsets of $\bbD$, 
\begin{equation}\lb{5.10}
\lim_{n\to\infty}\, \f{\psi_n^*(z)}{\varphi_n^*(z)} = F(z)
\end{equation}
\end{theorem} 

Clearly related to this is the following result of Golinskii-Nevai \cite{GN}: 

\begin{theorem} \lb{T5.3} Let $z\in\bbD$. Then 
\begin{equation}\lb{5.11}
\sum_{n=0}^\infty \, \biggl| \binom{\psi_n(z)}{-\psi_n^*(z)} + \beta  
\binom{\varphi_n(z)}{\varphi_n^*(z)}\biggr|^2 < \infty 
\end{equation}
if and only if 
\begin{equation}\lb{5.12}
\beta =F(z)
\end{equation} 
\end{theorem} 

From this point of view, $F$ is again the ``correct" analog of $m$! Indeed, the 
Golinskii-Nevai \cite{GN} proof uses Weyl limiting circles to prove the theorem 
(one is always in limit point case!). 

But this is not the end of the story. Define 
\begin{equation}\lb{5.13}
u_k =\psi_k +F(z)\varphi_k \qquad u_k^* = -\psi_k^* + F(z) \varphi_k^* 
\end{equation}
so $\binom{u_k}{u_k^*}$ is the unique solution of $\Xi_n =T_n(z) \Xi_0$ which 
is in $\ell^2$. In the OPRL case, the basic vector solution is of the form 
$\binom{u_n}{u_{n+1}}$, so we have the analog of \eqref{1.8e}, 
\begin{equation}\lb{5.14}
\ti m(z) = \f{u_0^*}{u_0} = \f{-1+F}{1+F} =zf
\end{equation}
So one analog of the $m$-function is $zf$. 

In particular, \eqref{5.14} implies 
\begin{equation}\lb{5.15}
\abs{u_k^*} < \abs{u_k}
\end{equation}
for $z\in\bbD$, and thus the rate of exponential decay of $\abs{\binom{u_k}{u_k^*}}$ is 
that of $u_k$. If there is such exponential decay in the sense that 
\begin{equation}\lb{5.16}
\gamma_2 = \lim_{n\to\infty}\, \biggl[\, \biggl\|\binom{u_n}{u_n^*}\biggr\|^{1/n} \, \biggr]
\end{equation}
exists, then, by \eqref{5.15}, 
\begin{equation}\lb{5.17}
\gamma_2 = \lim_{n\to\infty}\, \f{1}{n} \sum_{j=0}^{n-1} \log \abs{m_n^+} 
\end{equation}
where 
\begin{equation}\lb{5.18}
m_n^+ =\f{u_{n+1}}{u_n}
\end{equation}
For $n=0$, $u_1 =\psi_1 + F\varphi_1$, $u_0 =1+F$, $\psi_1 = \rho_0^{-1} (z+\bar\alpha_0)$, 
$\varphi_1 =\rho_0^{-1} (z-\bar\alpha_0)$, so by a direct calculation, 
\begin{equation}\lb{5.18a}
m_0^+ (z) = \rho_0^{-1} z(1-\bar\alpha_0 f)
\end{equation}
yet another reasonable choice for an $m$-function. 

Indeed, if $\gamma(z) =\lim_{n\to\infty} \f{1}{n}\log \|T_n(z)\|$ exists, the fact that 
$\det (T_n)=z^n$ implies that $\gamma =\log \abs{z}-\gamma_2$, and one finds in the case of 
stochastic Verblunsky coefficients that \cite{GTep,OPUC} 
\begin{equation}\lb{5.19}
\bbE (\log \abs{m_\omega^+(z)}) = \log \abs{z} -\gamma (z)
\end{equation}
an analog of a fundamental formula of Kotani \cite{Kot,S168} that in his case uses $m$! 

Finally, we turn to the connection of $m$ to whole-line Green's functions. Given $V$ on 
$(-\infty, \infty)$ and $z\in\bbC_+$, it is natural to look at the two solutions of 
\eqref{1.1}, $u_\pm (x,z)$, which are $\ell^2$ on $\pm (0,\infty)$ and the $m$-functions,  
\begin{equation}\lb{5.20}
m_\pm (z) = \pm \,\f{u'_\pm (0,z)}{u_\pm (0,z)}
\end{equation}
$m_\pm$ are the $m$-functions for $V(\pm x)\restriction [0,\infty)$. Standard Green's function 
formulae show that the integral kernel, $G(x,y;z)$ of $(-\f{d^2}{dx^2}+V-z)^{-1}$ is 
\[
G(x,y;z) =\f{u_-(x_<) u_+(x_>)}{(u_+(0) u'_-(0) - u'_+(0) u_-(0))}
\]
where $x_< =\min (x,y)$ and $x_> =\max (x,y)$. In particular, 
\begin{equation}\lb{5.21}
G(0,0;z) = -(m_+(z) + m_-(z))^{-1}
\end{equation}

A complete description of the OPUC analog would require too much space, so we sketch the ideas, 
leaving the details to \cite{OPUC}. Just as the difference equation is associated to a triagonal 
selfadjoint matrix whose spectral measure is the one generating the OPRL, any set of $\alpha$'s 
is associated to a five-diagonal unitary matrix, called the CMV matrix, whose spectral 
measure is the $d\mu$ with $\alpha_j (d\mu)=\alpha_j$. 

The CMV matrix is one-sided, but given $\{\alpha_j\}_{j=-\infty}^\infty$, one can define a 
two-sided CMV matrix, $\calE$, in a natural way. If $G(z)$ is the $00$ matrix element of 
$(\calE-z)^{-1}$, then (see \cite{GTep,Kh2000,OPUC}) 
\begin{equation}\lb{5.22}
G(z) = \f{f_+(z) f_-(z)}{1-zf_+(z) f_-(z)}
\end{equation}
where $f_+$ is the Schur function for $(\alpha_0, \alpha_1, \alpha_2, \dots)$ and $f_-$ the 
Schur function for $(-\bar\alpha_{-1}, -\bar\alpha_{-2}, \dots )$. On the basis of the analogy 
between \eqref{5.22} and \eqref{5.21}, Geronimo-Teplyaev \cite{GTep} called $f_+$ and 
$zf_-$ the $m_+$ and $m_-$ functions.

\medskip 
\section{Summary} \lb{s6} 

We have thus seen that there are many analogs of the $m$-function in the theory of OPUC: 
\begin{SL} 
\item[(1)] The Carath\'eodory function, $F(z)$, given by \eqref{3.1}, an analog of \eqref{1.2} 
and also related to the classic Weyl definition \eqref{5.11}/\eqref{5.12}.
\item[(2)] The Schur function, $f(z)$, given by \eqref{3.7} with a recursion, \eqref{3.8}, 
closer to the recursion \eqref{1.10} for the $m$-function of OPRL. $f$ also enters via 
\eqref{5.22}. 
\item[(3)] $zf(z)$, the $\ti m$-function of \eqref{5.14}. 
\item[(4)] The relative Szeg\H{o} function, \eqref{4.1}, which, via \eqref{4.3} and \eqref{1.11}, 
is an analog of $a_1 m(z)$. 
\item[(5)] The $m^+$-function, \eqref{5.18a}, which plays the role that $m$ does in Kotani theory. 
\end{SL}


\end{document}